\documentclass[10pt,twocolumn,twoside]{IEEEtran}

\usepackage[T1]{fontenc}
\usepackage{units}
\usepackage{amsthm}
\usepackage{amsmath}
\usepackage{amssymb}
\usepackage{enumerate,color,cite}
\usepackage{ragged2e}
\usepackage{graphicx}
\usepackage{algorithm, algorithmic}
\usepackage[hidelinks]{hyperref}

\graphicspath{{./Figs/}}

\makeatletter
\theoremstyle{plain}
\newtheorem{prop}{\propositionname}
\theoremstyle{plain}
\newtheorem{lem}{\lemmaname}
\theoremstyle{plain}
\newtheorem{thm}{\theoremname}
\theoremstyle{plain}
\newtheorem{cor}{\corollaryname}

\theoremstyle{definition}
\newtheorem{defn}{\definitionname}
\theoremstyle{remark}
\newtheorem{rem}{\remarkname}

\newtheorem{assumption}{Assumption}

\makeatother

\providecommand{\corollaryname}{Corollary}
\providecommand{\definitionname}{Definition}
\providecommand{\lemmaname}{Lemma}
\providecommand{\remarkname}{Remark}
\providecommand{\theoremname}{Theorem}
\providecommand{\propositionname}{Proposition}

\usepackage{caption}% http://ctan.org/pkg/caption
\newcommand*{\TitleFont}{%
      \fontsize{22}{26}%
      \selectfont}

\newcommand{\K}{\mathcal{K}}

\newcommand{\LL}{\mathcal{L}}
\newcommand{\KL}{\mathcal{KL}}

\begin{document}
\title{\TitleFont Generic Stability Implication from Full Information Estimation
to Moving-Horizon Estimation}
\author{Wuhua~Hu\thanks{W. Hu is with Towngas Energy, Shenzhen, China. The research was conducted
in his free time when he was working with Envision Digital in Singapore.
E-mail: hu.wuhua@towngas.com.cn}}
\maketitle
\begin{abstract}
Optimization-based state estimation is useful for handling of constrained
linear or nonlinear dynamical systems. It has an ideal form, known
as full information estimation (FIE) which uses all past measurements
to perform state estimation, and also a practical counterpart, known
as moving-horizon estimation (MHE) which uses most recent measurements
of a limited length to perform the estimation. Due to the theoretical
ideal, conditions for robust stability of FIE are relatively easier
to establish than those for MHE, and various sufficient conditions
have been developed in literature. This work reveals a generic link
from robust stability of FIE to that of MHE, showing that the former
implies at least a weaker robust stability of MHE which implements
a long enough horizon. The implication strengthens to strict robust
stability of MHE if the corresponding FIE satisfies a mild Lipschitz
continuity condition. The revealed implications are then applied to
derive new sufficient conditions for robust stability of MHE, which
further reveal an intrinsic relation between the existence of a robustly
stable FIE/MHE and the system being incrementally input/output-to-state
stable. 
\end{abstract}

\begin{IEEEkeywords}
Nonlinear systems; moving-horizon estimation; full information estimation;
state estimation; disturbances; robust stability; incremental input/output-to-state
stability 
\end{IEEEkeywords}

\thispagestyle{empty}

\section{Introduction}

Optimization-based state estimation is an estimation approach which
performs state estimation by solving an optimization problem. Compared
to conventional approaches like Kalman filtering (KF) which deals
with linear dynamical systems, and its extensions like extended KF
and unscented KF which can deal with nonlinear dynamical systems based
on linearization techniques, optimization-based approach has the advantage
of handling linear and nonlinear dynamical systems directly and also
including various physical or operational constraints \cite{rawlings2020model}.
The optimization formulation also admits flexible definition of the
objective function to be optimized, which in some cases is necessary
to accurately recover the state \cite{Hu2017robust}.

Optimization-based state estimation generally takes one of the two
forms: full information estimation (FIE) which uses all past measurements
to perform the estimation, and moving-horizon estimation (MHE) which
uses most recent measurements of a limited length to perform the estimation.
MHE is a practical approximate of FIE which is ideal and computationally
intractable. The interest in FIE lies in two aspects: serving as a
benchmark to MHE and providing useful insights for the stability analysis
of MHE. Recent studies on FIE and MHE are concentrated on analyses
of their stability and robustness, as are critical to guide the applications
of MHE. While earlier literature assumes restrictive and idealistic
conditions such as observability and/or zero or a priori known convergent
disturbances \cite{michalska1993moving,muske1993receding,michalska1995moving,rao2003constrained,alessandri2008moving,Rawlings2009Model,alessandri2010advances},
recent literature considers more practical conditions like detectability
and/or in the presence of bounded disturbances \cite{liu2013moving,rawlings2012optimization,rawlings2020model}.

An important progress was made in \cite{rawlings2012optimization},
which introduced the incremental input/output-to-state stability (i-IOSS)
concept on detectability of nonlinear systems (as developed in \cite{sontag1997output})
to study robust stability of FIE and MHE. Given an i-IOSS system,
it was shown that under mild conditions FIE is robustly stable and
convergent for convergent disturbances. Whereas, conditions for FIE
and MHE to be robustly stable under bounded disturbances were posted
as an open research challenge. Reference \cite{hu2015optimization}
provided a prompt response to this challenge,\footnote{The main results were obtained in the late of 2012 though the paper
was not able to be published until 2015.} identifying a general set of conditions for FIE to be robustly stable
for i-IOSS systems. Onwards, a series of researches have been inspired
to close the challenge.

A particular type of cost functions with a max-term was investigated
for a class of i-IOSS systems in \cite{ji2016robust}, establishing
robust stability of the FIE. The conditions were enhanced in \cite{muller2016nonlinear},
enabling robust stability also of MHE if a sufficiently long horizon
is applied. The conclusion was extended to MHE without a max-term
in the cost function \cite{muller2017nonlinear}. Meanwhile, it was
shown that MHE is convergent for convergent disturbances, with or
without a max-term. On the other hand, reference \cite{Hu2017robust}
revealed an implication link from robust stability of FIE to that
of MHE, and consequently identified rather general conditions for
MHE to be robustly stable by inheriting conditions which ensure robust
stability of the corresponding FIE. By making use of a Lipschitz continuity
condition introduced in \cite{Hu2017robust}, reference \cite{Allan2019inBook}
streamlined and generalized the analysis and results of \cite{muller2017nonlinear},
showing that the key is essentially to assume that the system satisfies
a form of exponential detectability. When global exponential detectability
is assumed, the MHE can further be shown to be robustly globally exponentially
stable by implementing properly time-discounted stage costs \cite{knufer2018robust}.

The reviewed robust stabilities of FIE and MHE were concluded based
on the concept defined in \cite{rawlings2012optimization}, which
is however found to be flawed in that such defined robustly stable
estimator does not imply convergence for convergent disturbances \cite{allan2021nonlinear,knufer2020time}.
This motivates a necessary modification of the stability concept in
\cite{allan2021nonlinear}, which redefines the estimate error bound
with a worst-case time-discounted instead of uniformly-weighted impact
of the disturbances. The new concept is an enhancement of the old
one, and was shown to be compatible with the original i-IOSS detectability
of a system, which is necessary for establishing robust stability
of both FIE and MHE. With this new concept, it becomes straightforward
to understand the earlier robust local/global exponential stability
results reported of FIE and MHE \cite{muller2017nonlinear,knufer2018robust}.
It also motivates some new results as developed in \cite[Ch. 4]{rawlings2020model},
which introduced a kind of stabilizability condition to establish
robust stability of FIE and MHE. Despite conceptual elegance, the
new condition involves an inequality which is uneasy to verify in
general.

Motivated by the new and stronger concept of robust stability, this
work aims to motivate general sufficient conditions for robust stability
of MHE by firstly establishing a generic implication link from robust
stability of FIE to that of MHE, and then transforming the challenge
into identifying sufficient conditions for ensuring robust stability
of FIE. While the reasoning approach is inspired by the ideas introduced
in \cite{Hu2017robust}, the contributions of this work are three-fold:
\begin{itemize}
\item A generic implication is established from robust stability of FIE
to practical robust stability of the corresponding MHE, and the implication
becomes stronger to robust stability of the MHE if the FIE admits
a certain Lipschitz continuity property;
\item Explicit computations of sufficient MHE horizons are provided for
the aforementioned implications to be effective in a local or global
sense of robust stability. Analyses on two exemplary cases further
illustrate that the upper bound of the MHE estimate error can eventually
decrease with respect to (w.r.t.) the horizon size once it is large
enough and satisfies certain conditions. 
\item Given the revealed implications, new sufficient conditions are derived
for robust stability of MHE by firstly developing those for the corresponding
FIE. An interesting finding is that a system being i-IOSS is necessary
but also sufficient for existence of a robustly stable FIE. Consequently
is a similar but weaker conclusion applicable to MHE.
\end{itemize}
The remaining of this work is organized as follows. Sec. \ref{sec: Notation-and-Preliminaries}
introduces notation, setup and necessary preliminaries. Sec. \ref{sec: optimization-based estimation}
defines general forms of FIE and MHE, and introduces robust stability
concepts. Sec. \ref{sec: stability link} reveals the implication
from global (or local) robust stability of FIE to that of the corresponding
MHE. Sec. \ref{sec: RGAS of FIE} applies the implication to establish
robust stability of MHE, by firstly developing conditions for ensuring
robust stability of the corresponding FIE. Finally, Sec. \ref{sec:Conclusion}
concludes the work.

\section{Notation, Setup and Preliminaries \label{sec: Notation-and-Preliminaries}}

The notation mostly follows the convention in \cite{Hu2017robust,rawlings2020model}.
The symbols $\mathbb{R}$, $\mathbb{R}_{\ge0}$ and $\mathbb{I}_{\ge0}$
denote the sets of real numbers, nonnegative real numbers and nonnegative
integers, respectively, and $\mathbb{I}_{a:b}$ denotes the set of
integers from $a$ to $b$. The constraints $t\ge0$ and $t\in\mathbb{I}_{\ge0}$
are used interchangeably to refer to the set of discrete times. The
symbol $\left|\cdot\right|$ denotes the Euclidean norm of a vector.
The bold symbol $\boldsymbol{x}_{a:b}$, denotes a sequence of vector-valued
variables $(x_{a},\,x_{a+1},\,...,\thinspace x_{b})$, and with a
function $f$ acting on a vector $x$, $f(\boldsymbol{x}_{a:b})$
stands for the sequence of function values $(f(x_{a}),\,f(x_{a+1}),\,...,\thinspace f(x_{b}))$.
Given a scalar function $f:\mathbb{R}\times\mathbb{R}\rightarrow\mathbb{R}$,
any $t\in\mathbb{R}$ and any scalar function $g$, the notation $f^{\circ n}(g(\cdot),t)$
refers to the $n$ fold composition of function $g$ by function $f$
subject to the given second argument $t$, and $f^{\circ0}$ equals
the identity function.

Throughout the paper, $t$ refers to a discrete time, and as a subscript
it indicates dependence on time $t$. Whereas, the subscripts or superscripts
$x$, $w$ and $v$ are used exclusively to indicate a function or
variable that is associated with the state ($x$), process disturbance
($w$) or measurement noise ($v$). The symbols $\left\lfloor x\right\rfloor $
and $\left\lceil x\right\rceil $ refer to the integers that are closest
to $x$ from below and above, respectively. Given two scalars $a$
and $b$, let $a\oplus b:=\max\{a,\,b\}$. The operator $\oplus$
is both associative and commutative, i.e., $(a\oplus b)\oplus c=a\oplus(b\oplus c)$
and $a\oplus b=b\oplus a$, and furthermore is distributive with respect
to (w.r.t.) increasing functions. That is, $\alpha(a\oplus b\oplus c)=\alpha(a)\oplus\alpha(b)\oplus\alpha(c)$
if the function $\alpha$ is increasing in the argument. The frequently
used $\K$, $\LL$ and $\KL$ functions are defined as follows. 
\begin{defn}
($\K$, $\LL$ and $\KL$ functions) A function $\alpha:\mathbb{R}_{\ge0}\to\mathbb{R}_{\ge0}$
is a $\K$ function if it is continuous, zero at zero, and strictly
increasing. A function $\varphi:\mathbb{R}_{\ge0}\to\mathbb{R}_{\ge0}$
is a $\LL$ function if it is continuous, nonincreasing and satisfies
$\varphi(t)\to0$ as $t\to\infty$. A function $\beta:\mathbb{R}_{\ge0}\times\mathbb{R}_{\ge0}\to\mathbb{R}_{\ge0}$
is a $\KL$ function if, for each $t\ge0$, $\beta(\cdot,t)$ is a
$\K$ function and for each $s\ge0$, $\beta(s,\cdot)$ is a $\LL$
function. 
\end{defn}
Consider a discrete-time system described by 
\begin{equation}
x_{t+1}=f(x_{t},w_{t}),\,\,y_{t}=h(x_{t})+v_{t},\label{eq:system}
\end{equation}
where $x_{t}\in\mathbb{X}\subseteq\mathbb{R}^{n}$ is the system state,
$w_{t}\in\mathbb{W}\subseteq\mathbb{R}^{g}$ the process disturbance,
$y_{t}\in\mathbb{Y}\subseteq\mathbb{R}^{p}$ the measurement, $v_{t}\in\mathbb{V}\subseteq\mathbb{R}^{p}$
the measurement disturbance, all at time $t$. Here we study state
estimation as an independent subject, and so control inputs (if there
were any) known up to the estimation time are treated as given constants,
which do not cause difficulty to later defined optimization and related
analyses and hence are neglected in the problem formulation for brevity
\cite{rawlings2012optimization,rawlings2020model}. The functions
$f$ and $h$ are assumed to be continuous and known. The initial
state $x_{0}$ and the disturbances $(w_{t},v_{t})$ are assumed to
be unknown but \textit{bounded}.

Consider two state trajectories. Let 
\begin{gather}
\begin{gathered}\pi_{0}:=x_{0}^{(1)}-x_{0}^{(2)},\\
\pi_{\tau+1}:=w_{\tau}^{(1)}-w_{\tau}^{(2)},\,\,\pi_{\tau+t+1}:=h(x_{\tau}^{(2)})-h(x_{\tau}^{(1)}),
\end{gathered}
\label{eq: pi}
\end{gather}
for all $\tau\in\mathbb{I}_{0:t-1}$, where $\pi_{\tau+t+1}$ corresponds
to $v_{\tau}^{(1)}-v_{\tau}^{(2)}$ if additive measurement noises
are present while identical measurements are assumed. Hence $\boldsymbol{\pi}_{0:2t}$
collects a sequence of deviation vectors, and its domain is denoted
as $\Pi$. Let $\iota(\pi_{\cdot})$ extract the time index (i.e.,
the original index $\tau$ above) of $\pi_{\cdot}$. Specifically,
we have 
\begin{equation}
\iota(\pi_{\tau+1})=\iota(\pi_{\tau+t+1})=\tau,\,\,\forall\tau\in\mathbb{I}_{-1:t-1}.\label{eq: tau}
\end{equation}
Note that the time index of $\pi_{0}$ is defined as $-1$ which refers
to the time associated with a priori information. 

With the above notation, we can have a concise statement of the i-IOSS
definition given in \cite{rawlings2020model,allan2021nonlinear}. 
\begin{defn}
\label{def: i-IOSS-1}(Concise definition of i-IOSS) The system $x_{t+1}=f(x_{t},w_{t})$,
$y_{t}=h(x_{t})$ is i-IOSS if there exists $\alpha\in\KL$ such that
for every two initial states $x_{0}^{(1)},x_{0}^{(2)}$, and two sequences
of disturbances $\boldsymbol{w}_{0:t-1}^{(1)},\boldsymbol{w}_{0:t-1}^{(2)}$,
the following inequality holds for all $t\in\mathbb{I}_{\ge0}$: 
\begin{equation}
\left|x_{t}^{(1)}-x_{t}^{(2)}\right|\le\max_{i\in\mathbb{I}_{0:2t}}\alpha\left(\left|\pi{}_{i}\right|,\,t-\iota(\pi_{i})-1\right),\label{eq: concise definition - i-IOSS}
\end{equation}
where $x_{t}^{(i)}$ is a shorthand of $x_{t}(x_{0}^{(i)},\boldsymbol{w}_{0:t-1}^{(i)})$
for $i\in\{1,\,2\}$. Furthermore, the system is exp-i-IOSS if in
the above inequality $\alpha$ admits an exponential form as $\alpha(s,\tau):=cs\lambda^{\tau}$
with certain $\lambda\in(0,1)$, $c>0$ and for all $s,\tau\ge0$. 
\end{defn}
As proved in \cite{allan2021nonlinear}, this definition of i-IOSS
is equivalent to the original one introduced in \cite{sontag1997output}
which applies $\K$ instead of $\KL$ function to each uncertainty
$|\pi_{i}|$ for all $i\ge1$.

The following defines Lipschitz continuity of a function at a given
point, and presents a related proposition. 
\begin{defn}
\label{def: Lipschitz-continuous} (Lipschitz continuity at a point)
Given a subset $\mathbb{S}\subseteq\mathbb{R}^{n}$, a function $f:\mathbb{R}^{n}\to\mathbb{R}^{m}$
is said to be Lipschitz continuous at a point $x^{*}\in\mathbb{S}$
over the subset $\mathbb{S}$ if there is a constant $c$ such that
$|f(x)-f(x^{*})|\le c|x-x^{*}|$ for all $x\in\mathbb{S}$. 
\end{defn}
\begin{prop}
\label{prop: T_lb}(Properties related to Lipschitz continuity) Given
any $\beta\in\KL$, $\eta\in(0,1)$, $\bar{s}>0$ and $\underline{s}\in(0,\,\bar{s}]$,
define $\tau_{\min}=\min\{\tau:\beta(s,\tau)\le\eta s,\,s\in[\underline{s},\bar{s}],\,\tau\in\mathbb{I}_{\ge0}\}$.
The following conclusions hold: a) the minimum $\tau_{\min}$ is well
defined, i.e., it exists; b) the inequality $\beta(s,\tau_{\min})\le\eta\cdot(s\oplus\underline{s})$
holds for all $s\in[0,\bar{s}]$; and c) if $\beta(\cdot,0)$ is Lipschitz
continuous at the origin over the range $[0,\,\bar{s}]$ , then a)
and b) hold also with $\underline{s}=0$. 
\end{prop}
\begin{IEEEproof}
a) With $\underline{s}>0$, if we take $\bar{\tau}:=\min\{\tau:\beta(\bar{s},\tau)\le\eta\underline{s}\}$,
which certainly exists and has a finite value by the definition of
a $\KL$ function, then we must have $\beta(s,\bar{\tau})\le\beta(\bar{s},\bar{\tau})\le\eta\underline{s}\le\eta s$
for all $s\in[\underline{s},\bar{s}]$. That is, $\bar{\tau}$ is
a feasible solution of $\tau$ and hence also an upper bound for the
considered minimization problem. Consequently, it is sufficient to
consider the decision variable $\tau$ in the finite set $\mathbb{I}_{0:\bar{\tau}}$.
By enumerating $\tau$ in the set $\mathbb{I}_{0:\bar{\tau}}$, the
inequality $\beta(s,\tau)\le\eta s$ can be checked for all $s\in[\underline{s},\bar{s}]$
subject to each given $\tau$. All valid values of $\tau$ can then
be collected in a set $\mathcal{T}$. Since $\mathcal{T}$ is nonempty
and finite, there must exist a minimum value, that is, $\tau_{\min}:=\min\{\tau:\tau\in\mathcal{T}\}$
is well defined. This completes the proof of a).

b) Since the inequality $\beta(s,\tau_{\min})\le\eta s=\eta\cdot(s\oplus\underline{s})$
holds for all $s\in[\underline{s},\bar{s}]$ by the definition of
$\tau_{\min}$, it is remaining to show that it holds also for all
$s\in[0,\underline{s})$. In the latter case, we have $\beta(s,\tau_{\min})<\beta(\underline{s},\tau_{\min})\le\eta\underline{s}=\eta\cdot(s\oplus\underline{s})$,
where the first inequality owes to the property of a $\KL$ function
and the second results from a). This completes the proof of b).

c) Since $\beta$ is a $\KL$ function, $\beta(\cdot,0)$ being Lipschitz
continuous at the origin over $[0,\bar{s}]$ implies that $\beta(s,\tau)\le a(\tau)s$
for all $s\in[0,\bar{s}]$ and any $\tau\in\mathbb{I}_{\ge0}$, where
$a(\tau)$ is the Lipschitz constant which is nonincreasing in $\tau$
and furthermore goes to zero as $\tau\rightarrow\infty$, i.e., $a\in\LL$.
Consequently, $\tau_{\min}=\min\{\tau:a(\tau)\le\eta,\,\tau\in\mathbb{I}_{\ge0}\}$
is well defined, which implies that $\beta(s,\tau_{\min})\le a(\tau_{\min})s\le\eta s$
for all $s\in[0,\bar{s}]$. Therefore, in this case both a) and b)
hold true for all $\underline{s}\ge0$, which includes the origin. 
\end{IEEEproof}

\section{Optimization-based State Estimation\label{sec: optimization-based estimation}}

Consider the system described in \eqref{eq:system}. Given a present
time $t$, the state estimation problem is to find an optimal estimate
of state $x_{t}$ based on historical measurements $\{y_{\tau}\}$
for $\tau$ in a time set. Ideally, all measurements up to time $t$
are used, leading to the so-called FIE; and practically, only measurements
are used within a limited distance backwards from time $t$, yielding
the so-called MHE. Both FIE and MHE can be cast as optimization problems.\footnote{Readers are referred to \cite{rawlings2012optimization} for a brief
introduction of their connection to control problems and difference
from probabilistic formulations.} To be concise, we will first define MHE and then treat FIE as a variant.

Let MHE implement a moving horizon of size $T$. The decision variables
are denoted as $(\boldsymbol{\chi}_{t-T:t},\boldsymbol{\omega}_{t-T:t-1},\boldsymbol{\nu}_{t-T:t-1})$,
which correspond to the system variables $(\boldsymbol{x}_{t-T:t},\boldsymbol{w}_{t-T:t-1},\boldsymbol{v}_{t-T:t-1})$.\footnote{As in \cite{rawlings2020model}, the last measurement is not considered
for ease of presentation, though the inclusion does not change the
conclusions.} And let the optimal decision variables be $(\hat{\boldsymbol{x}}_{t-T:t},\hat{\boldsymbol{w}}_{t-T:t-1},\hat{\boldsymbol{v}}_{t-T:t-1})$.
Since $\hat{\boldsymbol{x}}_{t-T+1:t}$ is uniquely determined from
$\hat{x}_{t-T}$ and $\hat{\boldsymbol{w}}_{t-T:t-1}$, the decision
variables essentially reduce to $(\chi_{t-T},\boldsymbol{\omega}_{t-T:t-1},\boldsymbol{\nu}_{t-T:t-1})$.

In addition, let $\bar{x}_{t-T}$ be a priori estimate of $x_{t-T}$,
and in particular $\bar{x}_{0}$ be bounded. Without loss of generality,
the prior estimates of the disturbances are assumed to be zero. Denote
the cost function as $V_{T}(\chi_{t-T}-\bar{x}_{t-T},\boldsymbol{\omega}_{t-T:t-1},\boldsymbol{\nu}_{t-T:t-1})$,
which penalizes uncertainties in the initial state, the process and
the measurements. Then, the MHE instance at time $t$ is defined by
the following optimization problem: 
\begin{equation}
\begin{aligned}\begin{array}{c}
\text{MHE (or FIE}\\
\text{if }T\leftarrow t):
\end{array} & \min V_{T}(\chi_{t-T}-\bar{x}_{t-T},\boldsymbol{\omega}_{t-T:t-1},\boldsymbol{\nu}_{t-T:t-1})\\
\text{s.t., } & \chi_{\tau+1}=f(\chi_{\tau},\omega_{\tau}),\,\,\forall\tau\in\mathbb{I}_{t-T:t-1},\\
 & y_{\tau}=h(\chi_{\tau})+\nu_{\tau},\,\,\forall\tau\in\mathbb{I}_{t-T:t-1},\\
 & \chi_{t-T}\in\mathbb{X},\,\boldsymbol{\omega}_{t-T:t-1}\in\mathbb{W}^{T},\,\boldsymbol{\nu}_{t-T:t-1}\in\mathbb{V}^{T}.
\end{aligned}
\label{eq: MHE}
\end{equation}
As $\boldsymbol{\nu}_{t-T:t-1}$ is uniquely determined by $\chi_{t-T}$
and $\boldsymbol{\omega}_{t-T:t-1}$, it is kept mainly for the convenience
of expressing the disturbance set and the objective function. Since
the global optimal solution $\hat{x}_{\tau}$, for any $\tau\le t$,
is dependent on time $t$ when the MHE instance is defined, to be
unambiguous we use $\hat{x}_{t}^{\star}$ to represent $\hat{x}_{t}$
that is solved from the instance defined at time $t$. This keeps
$\hat{x}_{t}^{\star}$ unchanged, while the realization $\hat{x}_{t}$
varies as MHE renews itself in time.

To define FIE, it suffices to adapt the horizon size $T$ in the MHE
formulation to taking the time-varying value $t$. The yielded FIE
has the form of an MHE but with complete data originating from the
zero initial time. To link them easily, an FIE is called the \textit{corresponding}
FIE of an MHE based on which the FIE is derived, and conversely is
the MHE called a \textit{corresponding} MHE of the FIE. Since it becomes
computationally intractable as time elapses, FIE is studied mainly
for its theoretical interest: its performance is viewed as a limit
or benchmark that MHE attempts to approach, and its stability can
be a good start point for analysis of MHE.

An important issue in designing FIE or MHE is to identify conditions
under which the associated optimization admits optimal estimates that
satisfy a robust stability property defined below. Let $\boldsymbol{x}_{0:t}(x_{0},\boldsymbol{w}_{0:t-1})$
denote a state sequence generated from an initial state $x_{0}$,
and a disturbance sequence $\boldsymbol{w}_{0:t-1}$. In addition,
define a bounded set for any given $\delta_{0}>0$: 
\begin{equation}
\mathbb{X}_{\delta_{0}}:=\{(x_{1},\,x_{2}):|x_{1}-x_{2}|\le\delta_{0},\,x_{1},x_{2}\in\mathbb{X}\}.\label{eq: local domain of x_0}
\end{equation}

\begin{defn}
\label{def: RGAS}(Robust stable estimation) The estimate $\hat{x}_{t}$
of state $x_{t}$ is based on partial or full sequence of the noisy
measurements, $\boldsymbol{y}_{0:t}=h(\boldsymbol{x}_{0:t}(x_{0},\boldsymbol{w}_{0:t-1}))+\boldsymbol{v}_{0:t}$.
The estimate is robustly asymptotically stable (RAS) if given any
$\delta_{0}>0$, there exist functions $\beta_{x},\,\beta_{w},\,\beta_{v}\in\KL$
such that the following inequality holds for all $(x_{0},\,\bar{x}_{0})\in\mathbb{X}_{\delta_{0}}$,
$\boldsymbol{\omega}_{0:t-1}\in\mathbb{W}^{t}$, $\boldsymbol{\nu}_{0:t-1}\in\mathbb{V}^{t}$
and $t\in\mathbb{I}_{\ge0}$: 
\begin{align}
\left|x_{t}-\hat{x}_{t}\right| & \le\beta_{x}(\left|x_{0}-\bar{x}_{0}\right|,t)\nonumber \\
 & \quad\oplus\max_{*\in\{w,v\}}\max_{\tau\in\mathbb{I}_{0:t-1}}\beta_{*}(\left|*_{\tau}\right|,\,t-\tau-1).\label{eq: RGAS condition}
\end{align}
The estimate is further said to be robustly globally asymptotically
stable (RGAS) if the inequality is satisfied for all $x_{0},\bar{x}_{0}\in\mathbb{X}$.
Moreover, if the estimate is RAS (or RGAS) and the $\KL$ function
admits an exponential form as $\beta_{*}(s,\tau):=c_{*}s\lambda^{t}$
with certain $\lambda\in(0,1)$, $c_{*}>0$ and for all $*\in\{x,w,v\}$,
then the estimate is said to be robustly exponentially stable (RES)
(or robustly globally exponentially stable (RGES)). 
\end{defn}
The last measurement $y_{t}$ and hence the corresponding noise $v_{t}$
is not considered in the above inequality, to keep the definition
consistent with the formulations of FIE and MHE. Here the definition
of RGAS strengthens the one introduced in \cite{rawlings2012optimization}
which applies $\K$ instead of $\KL$ functions to the disturbances.
This change is necessary to enable a desirable feature that a state
estimator which is RGAS must be convergent under convergent disturbances
\cite{rawlings2020model,allan2021nonlinear}. The next definition
presents a weaker alternative of the above robust stability, which
will also be needed in later analysis. 
\begin{defn}
\label{def: pRGAS}(Practical robust stable estimations) The estimate
defined in Definition \ref{def: RGAS} is practically RAS (pRAS) if given
any $\epsilon,\delta_{0}>0$, there exist functions $\beta_{x},\,\beta_{w},\,\beta_{v}\in\KL$
such that the following inequality holds for all $(x_{0},\,\bar{x}_{0})\in\mathbb{X}_{\delta_{0}}$,
$\boldsymbol{\omega}_{0:t-1}\in\mathbb{W}^{t}$, $\boldsymbol{\nu}_{0:t-1}\in\mathbb{V}^{t}$
and $t\in\mathbb{I}_{\ge0}$: 
\begin{align}
\left|x_{t}-\hat{x}_{t}\right| & \le\epsilon\oplus\beta_{x}(\left|x_{0}-\bar{x}_{0}\right|,t)\nonumber \\
 & \quad\oplus\max_{*\in\{w,v\}}\max_{\tau\in\mathbb{I}_{0:t-1}}\beta_{*}(\left|*_{\tau}\right|,\,t-\tau-1).\label{eq: RGAS condition-1}
\end{align}
The estimate is further said to be practically RGAS (pRGAS) if the inequality
is satisfied for all $x_{0},\bar{x}_{0}\in\mathbb{X}$. 
\end{defn}
The adverb ``practically'' before RGAS is employed to keep it in
line with the practical stability concept developed in control literature
(e.g., \cite{grune2017nonlinear}). Compared to an RGAS estimate,
the pRGAS estimate admits a looser bound with a non-vanishing constant
term $\epsilon$. As will be shown later, this term can be made arbitrarily
small if MHE implements a long enough horizon. In addition, it is
worthwhile to remark that there is no need to define such stability
variant for the case of RES or RGES because in either case the exponential
property ensures that the $\epsilon$ term will be fully tempered
and removed.

By applying the notation used in the concise definition of i-IOSS
in \eqref{eq: concise definition - i-IOSS}, conciser forms of the
above stability definitions can also be obtained, which will be useful
to convey stability conditions for both FIE and MHE later on. Towards
that, let $x_{0}^{(1)}:=x_{0}$, $w_{\tau}^{(1)}:=w_{\tau}$, $h(x_{\tau}^{(1)}):=y_{\tau}-v_{\tau}$,
$x_{0}^{(2)}:=\bar{x}_{0}$, $w_{\tau}^{(2)}:=0$, and $h(x_{\tau}^{(2)}):=y_{\tau}$
for all $\tau\in\mathbb{I}_{0:t-1}$. Consequently, the notation $\boldsymbol{\pi}_{0:2t}$
defined in \eqref{eq: pi} has a new realization, and its corresponding
domain is denoted as $\Pi$ if $\chi_{0},\,\bar{x}_{0}\in\mathbb{X}$
or $\Pi_{\delta_{0}}$ if $(\chi_{0},\,\bar{x}_{0})\in\mathbb{X}_{\delta_{0}}$.

With the new notation, it is straightforward to prove the following
equivalent definition of robust stabilities. 
\begin{lem}
\label{def: RGAS-1}(Concise equivalent definition of robust stable
estimations) The estimate defined in Definition \ref{def: RGAS} is
RAS if and only if given any $\delta_{0}>0$, there exists a function
$\beta\in\KL$ such that the following inequality holds for all $\boldsymbol{\pi}_{0:2t}\in\Pi_{\delta_{0}}$
and $t\in\mathbb{I}_{\ge0}$: 
\begin{equation}
\left|x_{t}-\hat{x}_{t}\right|\le\max_{i\in\mathbb{I}_{0:2t}}\beta(\left|\pi_{i}\right|,t-\iota(\pi_{i})-1).\label{eq: RGAS condition-1-1}
\end{equation}
And it is RGAS if and only if the above inequality holds for all $\boldsymbol{\pi}_{0:2t}\in\Pi$.
The estimate being RAS (or RGAS) is further said to be RES (or RGES),
if and only if the $\KL$ function $\beta$ admits the form as $\beta(s,\tau):=cs\lambda^{\tau}$
with certain $c>0$, $\lambda\in(0,1)$ and for all $\tau\ge0$ and
$s$ in the corresponding domain. 
\end{lem}
Similar concise definitions of pRAS and pRGAS can both be obtained
by adding a constant $\epsilon>0$ to the right hand side of \eqref{eq: RGAS condition-1-1},
in which the $\KL$ function $\beta$ will then be dependent on $\epsilon$.
For brevity, we do not state them with another lemma.

\section{Stability Implication from FIE to MHE\label{sec: stability link}}

At any discrete time, an MHE instance can be interpreted as the corresponding
FIE initiating from the start of the horizon over which the MHE instance
is defined. Thus, the corresponding FIE being robustly stable implies
that each MHE instance is robustly stable within the time horizon
over which the instance is defined. If we interpret this as MHE being
\textit{instance-wise} robustly stable, then the challenge reduces
to identifying conditions under which \textit{instance-wise} robust
stability implies robust stability of MHE. This observation was made
in \cite{Hu2017robust}, and the challenge was solved there for a
weaker definition of RGAS. This section resolves the challenge subject
to the stronger stability concept given by Definition \ref{def: RGAS}.

To that end, as in \cite{Hu2017robust}, we apply an ordinary assumption
on the prior estimate $\bar{x}_{t-T}$ of the initial state $x_{t-T}$
of an MHE instance.

\begin{assumption} \label{assump: A1} Given any time $t\ge T+1$,
the prior estimate $\bar{x}_{t-T}$ of $x_{t-T}$ is given such that
\[
|x_{t-T}-\bar{x}_{t-T}|\le|x_{t-T}-\hat{x}_{t-T}^{\star}|.
\]
\end{assumption}

The assumption is obviously true if $\bar{x}_{t-T}$ is set to $\hat{x}_{t-T}^{\star}$,
which is the MHE estimate obtained at time $t-T$. Alternatively,
a better $\bar{x}_{t-T}$ might be obtained with smoothing techniques
which use measurements both before and after time $t-T$ \cite{aravkin2016generalized,rawlings2020model}.

Next, we present an important lemma which links global robust stability
of MHE with that of its corresponding FIE. 
\begin{lem}
\label{lem: FIE-MHE-aRGAS and RGES}(Global stability implication
from FIE to MHE) Consider MHE under Assumption \ref{assump: A1}.
The following two conclusions hold:

a) (RGAS-->pRGAS/RGAS/RGES) If FIE is RGAS as per \eqref{eq: RGAS condition},
then there exists $\underline{T}\in\mathbb{I}_{\ge0}$ such that the
corresponding MHE under Assumption \ref{assump: A1} is pRGAS for
all $T\ge\text{\ensuremath{\underline{T}}}$. If further the $\K$
function $\beta_{x}(\cdot,0)$ in \eqref{eq: RGAS condition} is globally
Lipschitz continuous at the origin, then the implication strengthens
to the existence of $\underline{T}'\in\mathbb{I}_{\ge0}$ such that
the MHE is RGAS for all $T\ge\text{\ensuremath{\underline{T}}}'$;
and if furthermore the $\K$ functions $\beta_{w}(\cdot,0)$ and $\beta_{v}(\cdot,0)$
are also globally Lipschitz continuous at the origin, then the MHE
is RGES for all $T\ge\text{\ensuremath{\underline{T}}}'$.

b) (RGES-->RGES) If FIE is RGES as per \eqref{eq: RGAS condition},
then there exists $\underline{T}\in\mathbb{I}_{\ge0}$ such that the
corresponding MHE is RGES for all $T\ge\text{\ensuremath{\underline{T}}}$.
In particular, with the $\KL$ function $\beta_{x}(s,t):=c_{x}s\lambda^{t}$
in \eqref{eq: RGAS condition} for certain $c_{x}>0$ and $\lambda\in(0,1)$,
it is feasible to set $\underline{T}=1\oplus(\left\lfloor \log_{\lambda}\frac{1}{c_{x}}\right\rfloor +1)$. 
\end{lem}
\begin{IEEEproof}
a)\emph{ }\textit{RGAS-->pRGAS}\emph{.} Let $n:=\left\lfloor \frac{t}{T}\right\rfloor $,
and so $0\le t-nT\le T-1$. For all $\tau\in\mathbb{I}_{0:T-1}$,
the MHE and the corresponding FIE estimates are the same, both denoted
as $\hat{x}_{\tau}^{\star}$. So, given any $\tau\in\mathbb{I}_{0:t-nT}$,
the absolute estimation error $\left|x_{\tau}-\hat{x}_{\tau}^{\star}\right|$
satisfies the RGAS inequality given by \eqref{eq: RGAS condition}.
That is, we have 
\begin{align}
 & \left|x_{\tau}-\hat{x}_{\tau}^{\star}\right|\le\beta_{x}(|x_{0}-\bar{x}_{0}|,\thinspace\tau)\nonumber \\
 & \quad\oplus\max_{*\in\{w,v\}}\max_{\tau'\in\mathbb{I}_{0:\tau-1}}\beta_{*}(\left|*_{\tau'}\right|,\,\tau-\tau'-1),\label{eq: error bound from 0 to t-nT}
\end{align}
for all $\tau\in\mathbb{I}_{0:t-nT}$. Next, we proceed to show that
the RGAS property is maintained for all $\tau\in\mathbb{I}_{t-nT+1:t}$.

Let $t_{-1}:=0$ and $t_{i}:=t-(n-i)T$ for all $i\in\mathbb{I}_{0:n}$.
The MHE instance at time $t_{1}$ can be viewed as the corresponding
FIE confined to time interval $[t_{0},\thinspace t_{1}]$. Thus, the
MHE satisfies the RGAS property within this interval. That is, by
\eqref{eq: RGAS condition} we have: 
\begin{align*}
 & \left|x_{t_{1}}-\hat{x}_{t_{1}}^{\star}\right|\le\beta_{x}(\left|x_{t_{0}}-\hat{x}_{t_{0}}^{\star}\right|,\thinspace T)\\
 & \quad\oplus\max_{*\in\{w,v\}}\max_{\tau\in\mathbb{I}_{t_{0}:t_{1}-1}}\beta_{*}(\left|*_{\tau}\right|,\,t_{1}-\tau-1),
\end{align*}
where Assumption \ref{assump: A1} has been applied to produce the
first term of the right hand side of the inequality. Repeat this reasoning
for the MHE instance defined at time $t_{2}$ and then apply the above
inequality, yielding 
\begin{align*}
 & \left|x_{t_{2}}-\hat{x}_{t_{2}}^{\star}\right|\le\beta_{x}(\left|x_{t_{1}}-\hat{x}_{t_{1}}^{\star}\right|,\thinspace T)\\
 & \quad\oplus\max_{*\in\{w,v\}}\max_{\tau\in\mathbb{I}_{t_{1}:t_{2}-1}}\beta_{*}(\left|*_{\tau}\right|,\,t_{2}-\tau-1)\\
 & \le\beta_{x}^{\circ2}(\left|x_{t_{0}}-\hat{x}_{t_{0}}^{\star}\right|,\thinspace T)\\
 & \quad\oplus\max_{*\in\{w,v\}}\max_{\tau\in\mathbb{I}_{t_{0}:t_{1}-1}}\beta_{x}\left(\beta_{*}(\left|*_{\tau}\right|,\,t_{1}-\tau-1),\,T\right)\\
 & \quad\oplus\max_{*\in\{w,v\}}\max_{\tau\in\mathbb{I}_{t_{1}:t_{2}-1}}\beta_{*}(\left|*_{\tau}\right|,\,t_{2}-\tau-1).
\end{align*}

By induction, we obtain 
\begin{align*}
\left|x_{t}-\hat{x}_{t}^{\star}\right|= & \left|x_{t_{n}}-\hat{x}_{t_{n}}^{\star}\right|\le\beta_{x}^{\circ n}(\left|x_{t_{0}}-\hat{x}_{t_{0}}^{\star}\right|,\thinspace T)\\
 & \oplus\max_{*\in\{w,v\}}\max_{i\in\mathbb{I}_{0:n-1}}\max_{\tau\in\mathbb{I}_{t_{n-i-1}:t_{n-i}-1}}\\
 & \quad\beta_{x}^{\circ i}\left(\beta_{*}(\left|*_{\tau}\right|,\,t_{n-i}-\tau-1),\,T\right).
\end{align*}
Since $\left|x_{t_{0}}-\hat{x}_{t_{0}}^{\star}\right|$ satisfies
inequality \eqref{eq: error bound from 0 to t-nT}, this implies that
\begin{align}
\left|x_{t}-\hat{x}_{t}^{\star}\right|\le & \beta_{x}^{\circ n}\left(\beta_{x}(|x_{0}-\bar{x}_{0}|,\thinspace t_{0}),\,T\right)\nonumber \\
 & \oplus\max_{*\in\{w,v\}}\max_{i\in\mathbb{I}_{0:n}}\max_{\tau\in\mathbb{I}_{t_{n-i-1}:t_{n-i}-1}}\nonumber \\
 & \quad\beta_{x}^{\circ i}\left(\beta_{*}(\left|*_{\tau}\right|,\,t_{n-i}-\tau-1),\,T\right)\nonumber \\
= & \beta_{x}^{\circ\left\lfloor \frac{t}{T}\right\rfloor }\left(\beta_{x}(|x_{0}-\bar{x}_{0}|,\thinspace t-\left\lfloor \frac{t}{T}\right\rfloor T),\,T\right)\nonumber \\
 & \oplus\max_{*\in\{w,v\}}\max_{\tau\in\mathbb{I}_{0:t-1}}\nonumber \\
\beta_{x}^{\circ\left\lfloor \frac{t-\tau}{T}\right\rfloor } & \left(\beta_{*}(\left|*_{\tau}\right|,\,t-\left\lfloor \frac{t-\tau}{T}\right\rfloor T-\tau-1),\,T\right)\label{eq: key inequality}\\
=: & \beta_{x,T}(|x_{0}-\bar{x}_{0}|,\thinspace t)\nonumber \\
 & \oplus\max_{*\in\{w,v\}}\max_{\tau\in\mathbb{I}_{0:t-1}}\beta_{*,T}(\left|*_{\tau}\right|,\,t-\tau-1),\nonumber 
\end{align}
where the relation that $\left\lfloor \frac{t-\tau}{T}\right\rfloor =i$
for all $\tau\in\mathbb{I}_{t_{n-i-1}:t_{n-i}}$ and $i\in\mathbb{I}_{0:n}$
has been used to establish the first equality, and the three functions
$\beta_{x,T}$, $\beta_{w,T}$ and $\beta_{v,T}$ are defined by one-to-one
correspondence to the three preceding terms, each of which is dependent
on the horizon size $T$. By definition, it is obvious that $\beta_{x,T}(\cdot,\tau)$,
$\beta_{w,T}(\cdot,\tau)$ and $\beta_{v,T}(\cdot,\tau)$ are $\K$
functions for any given $\tau\ge0$. However, it is not necessary
that $\beta_{x,T}(s,\cdot)$, $\beta_{w,T}(s,\cdot)$ and $\beta_{v,T}(s,\cdot)$
are $\LL$ functions for all $s,T\ge0$, and so $\beta_{x,T}$, $\beta_{w,T}$
and $\beta_{v,T}$ are not necessarily $\KL$ functions. The remaining
proof is to show that there exists $\underline{T}$ such that $\beta_{x,T}$,
$\beta_{w,T}$ and $\beta_{v,T}$ are bounded by $\KL$ functions
plus a positive constant for all $T\ge\underline{T}$.

As $T$ increases from zero to infinity, the function instances form
a sequence $\boldsymbol{\beta}_{*,0:\infty}$, for each $*\in\{x,\,w,\,v\}$.
By the definition of $\beta_{*,T}$, given any $s$ and $\tau$, we
have $\lim_{T\rightarrow\infty}\beta_{*,T}(s,\tau)=\beta_{*}(s,\tau)$
for each $*\in\{x,\,w,\,v\}$. That is, $\beta_{*,T}$ point-wisely
converges to $\beta_{*}$ as $T\rightarrow\infty$. Next, we prove
that the convergence is uniform w.r.t. the arguments $s$ and $\tau$.
Since $|x_{0}-\bar{x}_{0}|$ is assumed to be bounded, there exists
$M>0$ such that $|x_{0}-\bar{x}_{0}|\le M$. Given any $0<\zeta<2\beta_{x}(M,0)$,
there must exist $T'$ satisfying $\beta_{x}(M\oplus\beta_{x}(M,\thinspace0),\thinspace T')<\frac{\zeta}{2}$,
such that for all $T\ge T'$ we will have for any $t<T$, 
\[
\sup_{x_{0},\bar{x}_{0},t}\left|\beta_{x,T}(|x_{0}-\bar{x}_{0}|,\thinspace t)-\beta_{x}(|x_{0}-\bar{x}_{0}|,\thinspace t)\right|=0<\zeta,
\]
and for any $t\ge T$, 
\begin{align*}
 & \sup_{x_{0},\bar{x}_{0},t}\left|\beta_{x,T}(|x_{0}-\bar{x}_{0}|,\thinspace t)-\beta_{x}(|x_{0}-\bar{x}_{0}|,\thinspace t)\right|\\
 & \le\beta_{x,T}(M,\thinspace t)+\beta_{x}(M,\thinspace t)<\beta_{x}^{\circ\left\lfloor \frac{t}{T}\right\rfloor }\left(\beta_{x}(M,\thinspace0),\,T\right)+\frac{\zeta}{2}\\
 & <\beta_{x}^{\circ\max\{0,\,\left\lfloor \frac{t}{T}\right\rfloor -1\}}\left(\frac{\zeta}{2},\,T\right)+\frac{\zeta}{2}\le\zeta.
\end{align*}
That is, $\sup_{x_{0},\bar{x}_{0},t}\left|\beta_{x,T}(|x_{0}-\bar{x}_{0}|,\thinspace t)-\beta_{x}(|x_{0}-\bar{x}_{0}|,\thinspace t)\right|<\zeta$
holds true for all $T\ge T'$. This implies that $\beta_{x,T}$ uniformly
converges to $\beta_{x}$. So can $\beta_{w,T}$ and $\beta_{v,T}$
be shown in the same way to uniformly converge to $\beta_{w}$ and
$\beta_{v}$, respectively.\footnote{If any of the uncertainties does not have a certain bound, then the
valid $T'$ will depend on the actual uncertainty magnitude and consequently
only semi-uniformly convergence will be proved.}

Consequently, given any $\epsilon>0$, there exists $\underline{T}>0$
such that $\beta_{x,T}(|x_{0}-\bar{x}_{0}|,\thinspace t)<\beta_{x}(|x_{0}-\bar{x}_{0}|,\thinspace t)+\epsilon/2$
and $\beta_{*,T}(\left|*_{\tau}\right|,\,t-\tau-1)<\beta_{*}(\left|*_{\tau}\right|,\,t-\tau-1)+\epsilon/2$
for all $*\in\{w,\,v\}$, $T\ge\underline{T}$, $x_{0},\bar{x}_{0}\in\mathbb{X}$,
$w_{\tau}\in\mathbb{W}$, $v_{\tau}\in\mathbb{V}$, $\tau\in\mathbb{I}_{0:t-1}$
and $t\in\mathbb{I}_{\ge0}$. Given each $*\in\{x,w,v\}$, let 
\[
\beta_{*}'(s,\tau):=\begin{cases}
\beta_{*}(s,\tau)+\epsilon/2 & \text{if }\beta_{*}(s,\tau)>\epsilon/2,\\
\beta_{*}(s,\tau) & \text{if }\beta_{*}(s,\tau)\le\epsilon/2,
\end{cases}
\]
for all $s,\tau\ge0$. Given any $*\in\{x,w,v\}$, it is easy to verify
that $\beta_{*}'$ is a $\KL$ function as $\beta_{*}$ is, and so
for any $T\ge\underline{T}$ we will have $\beta_{*,T}(s,\tau)\le\beta_{*}(s,\tau)+\epsilon/2\le\epsilon\oplus\beta_{*}'(s,\tau)$
for all $s,\tau$ in aforementioned domains. By applying these inequalities
to the last inequality of $\left|x_{t}-\hat{x}_{t}^{\star}\right|$
above, we conclude that given any $\epsilon>0$, there exists $\underline{T}>0$
such that 
\begin{align*}
\left|x_{t}-\hat{x}_{t}^{\star}\right| & \le\epsilon\oplus\beta_{x}'(|x_{0}-\bar{x}_{0}|,\thinspace t)\\
 & \quad\oplus\max_{*\in\{w,v\}}\max_{\tau\in\mathbb{I}_{0:t-1}}\begin{array}{c}
\beta_{*}'(\left|*_{\tau}\right|,\,t-\tau-1)\end{array}
\end{align*}
for all $T\ge\underline{T}$, where $\beta_{*}'\in\KL$ for all $*\in\{x,w,v\}$.
This implies that the MHE is pRGAS by definition, and hence completes
proof of the first conclusion of part a).

\textit{RGAS-->RGAS}\emph{.} If further $\beta_{x}(\cdot,0)$ from
\eqref{eq: RGAS condition} is globally Lipschitz continuous at the
origin, then for any $\tau\ge0$ there must exist a $\LL$ function
$\mu_{x}$ such that $\beta_{x}(s,\tau)\le\mu_{x}(\tau)s$ for all
$s\in\mathbb{X}$ because $\beta_{x}\in\mathcal{L}$. By applying
this property to \eqref{eq: key inequality}, for all $T\ge\text{\ensuremath{\underline{T}}}'$
satisfying $\mu_{x}(\text{\ensuremath{\underline{T}}}')\le\eta$ with
certain $\eta\in(0,1)$, the inequality there proceeds as 
\begin{align}
\left|x_{t}-\hat{x}_{t}^{\star}\right| & \le\mu_{x}^{\left\lfloor \frac{t}{T}\right\rfloor }(T)\mu(0)|x_{0}-\bar{x}_{0}|\oplus\max_{*\in\{w,v\}}\max_{\tau\in\mathbb{I}_{0:t-1}}\nonumber \\
 & \quad\mu_{x}^{\left\lfloor \frac{t-\tau}{T}\right\rfloor }(T)\beta_{*}(\left|*_{\tau}\right|,\,t-\left\lfloor \frac{t-\tau}{T}\right\rfloor T-\tau-1)\nonumber \\
 & \le\eta^{\frac{t}{T}-1}\mu_{x}(0)|x_{0}-\bar{x}_{0}|\nonumber \\
 & \quad\oplus\max_{*\in\{w,v\}}\max_{\tau\in\mathbb{I}_{0:t-1}}\eta^{\frac{t-\tau}{T}-1}\beta_{*}(\left|*_{\tau}\right|,\,0).\label{eq: RGAS->RGAS}
\end{align}
Consequently, the MHE is RGAS by definition, which completes the proof
of the second conclusion of part a).

\textit{RGAS-->RGES}\emph{.} If in addition to $\beta_{x}(\cdot,0)$,
the $\K$ functions $\beta_{w}(\cdot,0)$ and $\beta_{v}(\cdot,0)$
from \eqref{eq: RGAS condition} are globally Lipschitz continuous
at the origin, then inequality \eqref{eq: RGAS->RGAS} proceeds as
\begin{align*}
\left|x_{t}-\hat{x}_{t}^{\star}\right| & \le\eta^{\frac{t}{T}-1}\mu_{x}(0)|x_{0}-\bar{x}_{0}|\\
 & \quad\oplus\max_{*\in\{w,v\}}\max_{\tau\in\mathbb{I}_{0:t-1}}\eta^{\frac{t-\tau}{T}-1}\mu_{*,0}\left|*_{\tau}\right|,
\end{align*}
where $\mu_{*,0}$ is the Lipschitz constant of $\beta_{*}(\cdot,0)$
at the origin for all $*\in\{w,v\}$. Consequently, the MHE is RGES
by definition, which completes the third conclusion of part a). 

b)\emph{ }\textit{RGES-->RGES}\emph{.} In this case, we have $\beta_{*}(s,\tau)=c_{*}s\lambda^{\tau}$
with certain $c_{*}>0$ and $\lambda\in(0,1)$ for each $*\in\{x,w,v\}$.
Since RGES implies RGAS, the deduction in the proof of part a) is
applicable. In particular, the three functions $\beta_{x,T}$, $\beta_{w,T}$
and $\beta_{v,T}$ are derived explicitly as: 
\begin{equation}
\begin{aligned}\beta_{x,T}(|x_{0}-\bar{x}_{0}|,\thinspace t) & =c_{x}^{\left\lfloor \frac{t}{T}\right\rfloor +1}|x_{0}-\bar{x}_{0}|\lambda^{t},\\
\beta_{w,T}(\left|w_{\tau}\right|,\,t-\tau-1) & =c_{x}^{\left\lfloor \frac{t-\tau}{T}\right\rfloor }c_{w}\left|w_{\tau}\right|\lambda^{t-\tau-1},\\
\beta_{v,T}(\left|v_{\tau}\right|,\,t-\tau-1) & =c_{x}^{\left\lfloor \frac{t-\tau}{T}\right\rfloor }c_{v}\left|v_{\tau}\right|\lambda^{t-\tau-1},
\end{aligned}
\label{eq: RGES error bound funcs}
\end{equation}
for all $x_{0},\bar{x}_{0}\in\mathbb{X}$, $w_{\tau}\in\mathbb{W}$,
$v_{\tau}\in\mathbb{V}$, $\tau\in\mathbb{I}_{0:t-1}$ and $t\in\mathbb{I}_{\ge0}$.
If $c_{x}\le1$, it is obvious that $\beta_{x,T}$, $\beta_{w,T}$
and $\beta_{v,T}$ are all upper bounded by $\KL$ functions in exponential
forms, which implies RGES of the MHE for all $T\ge\underline{T}:=1\ge1+\left\lfloor \log_{\lambda}\frac{1}{c_{x}}\right\rfloor $.
And if $c_{x}>1$, by defining $\underline{T}=1+\left\lfloor \log_{\lambda}\frac{1}{c_{x}}\right\rfloor $
($\ge1$) we have $c_{x}^{\left\lfloor \frac{t}{T}\right\rfloor }\lambda^{t}\le c_{x}^{\frac{t}{T}}\lambda^{t}\le\left(c_{x}^{\frac{1}{\underline{T}}}\lambda\right)^{t}\le\left(\left(c_{x}\right)^{\frac{1}{\log_{\lambda}\frac{1}{c_{x}}}}\lambda\right)^{t}=1$
for all $T\ge\underline{T}$, where the penultimate equality holds
only for $t=0$. Consequently $\beta_{x,T}$, $\beta_{w,T}$ and $\beta_{v,T}$
are again all upper bounded by $\KL$ functions in exponential forms
which imply RGES of the MHE. Combining the two cases, we conclude
that the MHE is RGES for all $T\ge\underline{T}:=1\oplus(\left\lfloor \log_{\lambda}\frac{1}{c_{x}}\right\rfloor +1)$,
which completes the proof of b). 
\end{IEEEproof}
Lemma \ref{lem: FIE-MHE-aRGAS and RGES}.a) indicates that RGAS of
FIE implies pRGAS of the corresponding MHE which implements a long
enough horizon, and that the implication strengthens to RGAS or RGES
of the MHE if the FIE additionally satisfies certain global Lipschitz
continuity conditions. In the latter case, it is not certain that
the MHE will be RGES unless the three bound functions \{$\beta_{x}$,
$\beta_{w}$, $\beta_{v}$\} of the FIE all satisfy the global Lipschitz
conditions. Also note that the main conclusion of Lemma \ref{lem: FIE-MHE-aRGAS and RGES}.b)
is implied by the third conclusion of Lemma \ref{lem: FIE-MHE-aRGAS and RGES}.a)
and that its proof above is presented mainly for derivation of an
explicit bound of the sufficient horizon size.
\begin{rem}
\label{rmk: Lipschitz-exponential-FIE}In the second case of Lemma
\ref{lem: FIE-MHE-aRGAS and RGES}, an FIE being RGAS and its bound
function $\beta_{x}(\cdot,0)$ being globally Lipschitz at the origin
imply that the system is exp-i-IOSS. This can be proved by extending
the proof for the local case in \cite[Proposition 2]{Allan2019inBook}
to the global case. However, the two conditions do not necessarily
imply that the FIE will be RGES (or at least the proof is unclear
yet), though practical RGES is guaranteed by Lemma \ref{lem: FIE-MHE-aRGAS and RGES}
if $\beta_{w}(\cdot,0)$ and $\beta_{v}(\cdot,0)$ are also globally
Lipschitz at the origin. On the other hand, if an FIE is RGES, then
the Lipschitz conditions in the third case are valid and hence the
corresponding MHE will be RGES by Lemma \ref{lem: FIE-MHE-aRGAS and RGES},
as is in line with the existing results reported in \cite{knufer2018robust,allan2020lyapunov}.
\end{rem}
\begin{rem}
\label{rmk: relation-to-recent-work}Reference \cite{knuefer2021nonlinear}
proved existence of a finite horizon for MHE to be RGAS if the corresponding
FIE induces a global contraction for the estimation error which essentially
requires the dynamical system to be exp-i-IOSS. As referred to Remark
\ref{lem: FIE-MHE-aRGAS and RGES}, this conveys the same necessary
condition as in the second case of Lemma \ref{lem: FIE-MHE-aRGAS and RGES}.
Overall, Lemma \ref{lem: FIE-MHE-aRGAS and RGES} covers more general
conclusions, which links RGAS of FIE to pRGAS of the corresponding
MHE but also RGAS and RGES of MHE as two naturally induced cases subject
to extra mild conditions. 
\end{rem}
\begin{rem}
\label{rmk: relation-to-MPC-literature}Reference \cite{jadbabaie2005stability}
proved existence of a finite horizon for a model predictive control
scheme to exponentially stabilize an unperturbed nonlinear system
via state feedback. Though the result could be relevant to the strongest
conclusion of Lemma \ref{lem: FIE-MHE-aRGAS and RGES} when the FIE
is RGES, a strict and meaningful connection is yet to be explored
for establishing robustly stable estimators here when both disturbances
and i-IOSS stability are in place which are yet absent in the system
setup of \cite{jadbabaie2005stability}.
\end{rem}
Next, we present additional results when the conditions of Lemma \ref{lem: FIE-MHE-aRGAS and RGES}
are relaxed and confined to local regions of the estimation problem.
In this case, as a byproduct, an explicit bound for the sufficient
MHE horizon size is derived when the corresponding FIE is known to
be RAS or RGAS. 
\begin{lem}
\label{lem: FIE-MHE-RAS-aRAS}(Local/global stability implication
from FIE to MHE) Consider MHE under Assumption \ref{assump: A1}.
Let the disturbances be bounded as $|w_{t}|\le\delta_{w}$ and $|v_{t}|\le\delta_{v}$
for all $t\in\mathbb{I}_{\ge0}$. Given any $\eta\in(0,\thinspace1)$,
the following four conclusions hold:

a) (RAS-->pRAS/RAS) Given any $\delta_{0}>0$ and $\epsilon\in(0,\eta\bar{s}]$,
where $\bar{s}:=\beta_{x}(\delta_{0},\thinspace0)\oplus\beta_{w}(\delta_{w},\,0)\oplus\beta_{v}(\delta_{v},\,0)$,
if FIE is RAS as per \eqref{eq: RGAS condition}, then the corresponding
MHE is pRAS for all $(x_{0},\,\bar{x}_{0})\in\mathbb{X}_{\delta_{0}}$
and $T\ge\underline{T}(\epsilon):=\min\{\tau:\beta_{x}(s,\tau)\le\eta s,\,s\in[\frac{\epsilon}{\eta},\bar{s}]\}$.
If in addition $\beta_{x}(\cdot,\tau)$ is Lipschitz continuous at
the origin over $[0,\bar{s}]$ for all $\tau\ge0$, then the corresponding
MHE is RAS for all $(x_{0},\,\bar{x}_{0})\in\mathbb{X}_{\delta_{0}}$
and $T\ge\underline{T}(0)$.

b) (RGAS-->pRGAS/RGAS) If the conditions in a) are satisfied for
all $x_{0},\bar{x}_{0}\in\mathbb{X}$ (which replaces $(x_{0},\,\bar{x}_{0})\in\mathbb{X}_{\delta_{0}}$),
then the two conclusions in a) hold globally.

c) (RES-->RES) Given $\bar{s}$ defined in a) and any $\delta_{0}>0$,
if FIE is RES as per \eqref{eq: RGAS condition}, in which for each
$*\in\{x,w,v\}$ the $\KL$ function $\beta_{*}$ has an exponential
form as $\beta_{*}(s,\tau):=c_{*}s\lambda^{\tau}$ for certain $c_{*}>0$,
$\lambda\in(0,1)$ and all $s\in[0,\bar{s}]$ and $\tau\ge0$, then
the corresponding MHE is RES for all $(x_{0},\,\bar{x}_{0})\in\mathbb{X}_{\delta_{0}}$
and $T\ge\underline{T}:=1\oplus(\left\lfloor \log_{\lambda}\frac{1}{c_{x}}\right\rfloor +1)$.

d) (RGES-->RGES) If the conditions in c) are satisfied for all $x_{0},\bar{x}_{0}\in\mathbb{X}$,
then the conclusion in c) holds globally. 
\end{lem}
\begin{IEEEproof}
a) \textit{RAS-->pRAS.} Given any $\eta\in(0,\,1)$ and $\epsilon\in(0,\eta\bar{s}]$,
define $\tilde{\beta}_{x}(s)=\eta\cdot(s\oplus\frac{\epsilon}{\eta})$
for all $s\in[0,\bar{s}]$. The proof follows the same procedure of
the proof for Lemma \ref{lem: FIE-MHE-aRGAS and RGES}.a), but further
bounds the term after the function composition with $\beta_{x}\circ$
by $\tilde{\beta}_{x}(s)$ defined here in each induction step. That
the value of $\tilde{\beta}_{x}(s)$ lies in $[0,\bar{s}]$ makes
the same compose-and-then-bound approach applicable throughout the
induction steps.

Let $\underline{T}(\epsilon):=\min\{\tau:\beta_{x}(s,\tau)\le\eta s,\,s\in[\frac{\epsilon}{\eta},\bar{s}]\}$,
which is well defined by Proposition \ref{prop: T_lb}.a). For any
$T\ge\underline{T}(\epsilon)$, Proposition \ref{prop: T_lb}.b) is
applicable, and so the aforementioned induction leads to the following
bound of the MHE estimate error:{\small{}{} 
\begin{align}
 & \left|x_{t}-\hat{x}_{t}^{\star}\right|\le\epsilon\oplus\eta^{\left\lfloor \frac{t}{T}\right\rfloor }\beta_{x}(|x_{0}-\bar{x}_{0}|,\thinspace t-\left\lfloor \frac{t}{T}\right\rfloor T)\nonumber \\
 & \quad\oplus\max_{*\in\{w,v\}}\max_{\tau\in\mathbb{I}_{0:t-1}}\eta^{\left\lfloor \frac{t-\tau}{T}\right\rfloor }\beta_{*}(\left|*_{\tau}\right|,\,t-\left\lfloor \frac{t-\tau}{T}\right\rfloor T-\tau-1)\nonumber \\
 & \le\epsilon\oplus\frac{1}{\eta}\beta_{x}(|x_{0}-\bar{x}_{0}|,\thinspace0)\eta^{\frac{t}{T}}\nonumber \\
 & \quad\oplus\max_{*\in\{w,v\}}\max_{\tau\in\mathbb{I}_{0:t-1}}\frac{1}{\eta}\beta_{*}(\left|*_{\tau}\right|,\,0)\eta^{\frac{t-\tau}{T}}\nonumber \\
 & =:\epsilon\oplus\frac{1}{\eta}\beta_{x,T}(|x_{0}-\bar{x}_{0}|,\thinspace t)\nonumber \\
 & \quad\oplus\max_{*\in\{w,v\}}\max_{\tau\in\mathbb{I}_{0:t-1}}\frac{1}{\eta}\beta_{*,T}(\left|*_{\tau}\right|,\,t-\tau-1),\label{eq: error bound}
\end{align}
}where the three functions $\beta_{x,T}$, $\beta_{w,T}$ and $\beta_{v,T}$
are defined by one-to-one correspondence to the three preceding terms,
each of which is dependent on the horizon size $T$. With $\eta\in(0,1)$,
it is obvious that $\beta_{x,T}$, $\beta_{w,T}$ and $\beta_{v,T}$
are all $\KL$ functions. Consequently, the above estimate error bound
implies that the MHE is pRAS by definition.

\textit{RAS-->RAS.} When $\beta_{x}(\cdot,\tau)$ is Lipschitz continuous
at the origin over the domain $[0,\bar{s}]$ for all $\tau\ge0$,
by Proposition \ref{prop: T_lb}.c) the value $\underline{T}(\epsilon)$
is well defined also at $\epsilon=0$. As a consequence, the same
proof above will show that the implied pRAS collapses to RAS, as the
ambiguity caused by $\epsilon$ disappears.

b) \textit{RGAS-->pRGAS/RGAS.} The proof is straightforward by extending
the reasoning in the proof of a) to the entire applicable domain with
$x_{0},\bar{x}_{0}\in\mathbb{X}$.

c) \textit{RES-->RES.} In this case, we have $\beta_{x}(s,\tau)=c_{x}s\lambda^{\tau}$
with certain $c_{x}>0$, $\lambda\in(0,1)$ and for all $s\in[0,\bar{s}]$
and $\tau\ge0$. If $c_{x}\le1$, let $\underline{T}:=1$ and hence
$\beta_{x}(s,T)\le\beta_{x}(s,\text{\ensuremath{\underline{T}}})=c_{x}s\lambda<s$
for all $T\ge\underline{T}$ and $s\in[0,\bar{s}]$. And if $c_{x}>1$,
let $\text{\ensuremath{\underline{T}}}:=\left\lfloor \log_{\lambda}\frac{1}{c_{x}}\right\rfloor +1$
and hence $\beta_{x}(s,T)\le\beta_{x}(s,\text{\ensuremath{\underline{T}}})=c_{x}s\lambda^{\left\lfloor \log_{\lambda}\frac{1}{c_{x}}\right\rfloor +1}<c_{x}s\lambda^{\log_{\lambda}\frac{1}{c_{x}}}=s$
for all $T\ge\underline{T}$ and $s\in[0,\bar{s}]$. In either case,
the induction steps in the proof of Lemma \ref{lem: FIE-MHE-aRGAS and RGES}.b)
are applicable as the domain of the function composition via $\beta_{x}\circ$
is kept within $[0,\bar{s}]$ due to the contraction induced by the
inequality $\beta_{x}(s,T)<s$. As a consequence, the state estimate
error will again be bounded as per \eqref{eq: RGES error bound funcs}.
By applying the same argument right after \eqref{eq: RGES error bound funcs}
subject to $\underline{T}=1\oplus(\left\lfloor \log_{\lambda}\frac{1}{c_{x}}\right\rfloor +1)$
(which consolidates the two cases), we conclude that the MHE corresponding
to the FIE is RES. This completes the proof of c).

d) \textit{RGES-->RGES.} The proof is straightforward by extending
the reasoning in the proof of c) to the entire problem domain with
$x_{0},\bar{x}_{0}\in\mathbb{X}$. 
\end{IEEEproof}
\begin{rem}
Two remarks follow on the proof of Lemma \ref{lem: FIE-MHE-RAS-aRAS}.a).
a) The proof shows that the upper bound of the estimate error is scaled
by the factor $\eta^{\left\lfloor \frac{t}{T}\right\rfloor }$, where
$\eta$ is a value that bonds with the horizon bound $\underline{T}(\epsilon)$.
Given any $T\ge\underline{T}(\epsilon)$, the error bound can in fact
be tightened by replacing the factor $\eta^{\left\lfloor \frac{t}{T}\right\rfloor }$
with $\left(\eta\phi(T-\underline{T})\right)^{\left\lfloor \frac{t}{T}\right\rfloor }$,
where $\phi$ is a certain $\LL$ function (and $\underline{T}$ is
short for $\underline{T}(\epsilon)$). As shown in Appendix A, as
time $t$ is large enough, the derivative of this tighter factor will
be negative and hence the error bound will be decreasing w.r.t. $T$
if the horizon size $T$ and the associated bound $\underline{T}(\epsilon)$
satisfy certain conditions that depend on the $\KL$ function $\beta_{x}$.
b) The proof indicates that it is valid for MHE to change its horizon
size over time according to the amplitude of uncertainty in an intermediate
initial state, as long as sufficient contraction of the uncertainty
is maintained. This leaves space for adapting the horizon size of
MHE such that a flexible and time-varying tradeoff can be made between
estimation quality and computational complexity. 
\end{rem}
Lemma \ref{lem: FIE-MHE-RAS-aRAS} indicates that there always exists
a large enough horizon size $T$ such that RAS/RGAS of its corresponding
FIE implies pRAS/pRGAS of MHE, and further RAS/RGAS of the MHE if
the FIE admits a certain Lipschitz continuity property. When the corresponding
FIE is RES/RGES, the implication directly strengthens to RES/RGES
of the MHE.

Note that the implication for pRGAS of MHE in Lemma \ref{lem: FIE-MHE-RAS-aRAS}.b)
is equivalent to the one concluded in Lemma \ref{lem: FIE-MHE-aRGAS and RGES}.a),
but here it is derived from a different proof. The advantage is to
have an explicit way of computing a sufficient horizon size. The implication
for RGES of MHE in Lemma \ref{lem: FIE-MHE-RAS-aRAS}.d) is again
equivalent to that presented in Lemma \ref{lem: FIE-MHE-aRGAS and RGES}.b),
while here it is extended naturally from a local implication.

The bound $\underline{T}$ for a sufficient horizon size as given
in Lemma \ref{lem: FIE-MHE-RAS-aRAS}.a) can be derived explicitly
if the $\KL$ function $\beta_{x}$ has particular forms. Given $a_{1},a_{2}\ge1,\,b_{1}\in(0,1),\,b_{2}>0$
and $c_{1},c_{2}>0$, two expressions of $\underline{T}$ are derived
as follows: 
\[
\underline{T}=\begin{cases}
\left\lceil \log_{b_{1}}\frac{\eta}{c_{1}\bar{s}^{a_{1}-1}}\right\rceil , & \text{if }\beta_{x}(s,\tau):=c_{1}s^{a_{1}}b_{1}^{\tau},\\
\left\lceil \sqrt[b_{2}]{\frac{c_{2}\bar{s}^{a_{2}-1}}{\eta}}-1\right\rceil , & \text{if }\beta_{x}(s,\tau):=c_{2}s^{a_{2}}(\tau+1)^{-b_{2}}.
\end{cases}
\]
The given constants $\{a_{i},b_{i},c_{i}\}$ for $i\in\{1,2\}$ are
such that the two $\beta_{x}(\cdot,\tau)$'s above are Lipschitz continuous
at the origin over $[0,\,\bar{s}]$ for all $\tau\in\mathbb{I}_{\ge0}$.

\section{Robust Stability of FIE and MHE\label{sec: RGAS of FIE}}

This section presents new sufficient conditions for robust stability
of FIE, and also for that of the corresponding MHE by applying the
stability implication revealed in last section.

To make the conditions easy to interpret and the proof concise to
present, the following notations are introduced in the spirit of \eqref{eq: pi}:
\begin{equation}
\begin{aligned}\pi_{0} & :=x_{0}-\bar{x}_{0},\,\pi_{\tau+1}:=w_{\tau},\,\pi_{\tau+t+1}:=v_{\tau},\\
\hat{\pi}_{0} & :=x_{0}-\hat{x}_{0},\,\hat{\pi}_{\tau+1}:=w_{\tau}-\hat{w}_{\tau},\,\hat{\pi}_{\tau+t+1}:=v_{\tau}-\hat{v}_{\tau}\\
\tilde{\pi}_{0} & :=\chi_{0}-\bar{x}_{0},\,\tilde{\pi}_{\tau+1}:=\omega_{\tau},\,\tilde{\pi}_{\tau+t+1}:=\nu_{\tau},\\
\hat{\tilde{\pi}}_{0} & :=\hat{x}_{0}-\bar{x}_{0},\,\hat{\tilde{\pi}}_{\tau+1}:=\hat{w}_{\tau},\,\hat{\tilde{\pi}}_{\tau+t+1}:=\hat{v}_{\tau},
\end{aligned}
\label{eq: pi's}
\end{equation}
for all $\tau\in\mathbb{I}_{0:t-1}$, where $\hat{\pi}_{\cdot}$ and
$\hat{\tilde{\pi}}_{\cdot}$ refer to the optimal estimates of $\pi_{\cdot}$
and $\tilde{\pi}_{\cdot}$, respectively. Given any $*\in\{\pi,\,\hat{\pi},\,\tilde{\pi},\,\hat{\tilde{\pi}}\}$,
notation $\boldsymbol{*}_{0:2t}$ collects the sequence of vector
variables $(*_{\cdot})_{0:2t}$ and the corresponding domain is denoted
by $\Pi$ if $\chi_{0},\,\bar{x}_{0}\in\mathbb{X}$ and as $\Pi_{\delta_{0}}$
if $(\chi_{0},\,\bar{x}_{0})\in\mathbb{X}_{\delta_{0}}$, where $\mathbb{X}_{\delta_{0}}$
is defined in \eqref{eq: local domain of x_0}. Function $\iota(*_{\cdot})$
extracts the time index (i.e., the original index $\tau$) of $*_{\cdot}$
as per \eqref{eq: tau}. Given any $i\in\mathbb{I}_{0:2t}$ , it is
easy to verify that $\iota(\pi_{i})=\iota(\hat{\pi}_{i})=\iota(\tilde{\pi}_{i})=\iota(\hat{\tilde{\pi}}_{i})$
and that $\pi_{i}=\hat{\pi}_{i}+\hat{\tilde{\pi}}_{i}$. 

Given any $\delta_{0}>0$, the next two assumptions are introduced
to establish robust stability of FIE.

\begin{assumption} \label{assump: A2} There exists $\underbar{\ensuremath{\rho}},\rho\in\KL$
such that the cost function of FIE $V_{t}(\tilde{\boldsymbol{\pi}}_{0:2t})$,
which is continuous, satisfies the following inequality for all $\tilde{\boldsymbol{\pi}}_{0:2t}\in\Pi_{\delta_{0}}$
and $t\in\mathbb{I}_{\ge0}$: 
\begin{align}
 & \max_{i\in\mathbb{I}_{0:2t}}\underbar{\ensuremath{\rho}}(\left|\tilde{\pi}_{i}\right|,\,t-\iota(\tilde{\pi}_{i})-1)\le V_{t}(\tilde{\boldsymbol{\pi}}_{0:2t})\label{eq: Assumption 2}\\
 & \quad\le\max_{i\in\mathbb{I}_{0:2t}}\ensuremath{\rho}(\left|\tilde{\pi}_{i}\right|,\,t-\iota(\tilde{\pi}_{i})-1).\nonumber 
\end{align}
\end{assumption}

\begin{assumption} \label{assump: A3} $\KL$ function $\alpha$
from the i-IOSS property \eqref{eq: concise definition - i-IOSS}
and $\KL$ functions $\text{\ensuremath{\underbar{\ensuremath{\rho}}}}$
and $\rho$ from Assumption \ref{assump: A2} satisfy the following
inequality for all $\tilde{\boldsymbol{\pi}}_{0:2t}\in\Pi_{\delta_{0}}$
and $\tau,\,\tau'\in\mathbb{I}_{0:t}$: 
\begin{equation}
\alpha\left(2\underbar{\ensuremath{\rho}}^{-1}\left(\rho(|\tilde{\pi}_{i}|,\,\tau),\,\tau'\right),\,\tau'\right)\le\bar{\alpha}(|\tilde{\pi}_{i}|,\,\tau),\label{eq: Assumption 3}
\end{equation}
for certain $\bar{\alpha}\in\KL$, in which $\underbar{\ensuremath{\rho}}^{-1}(\cdot,\,\tau)$
is the inverse of $\underbar{\ensuremath{\rho}}(\cdot,\,\tau)$ w.r.t.
its first argument given the second argument $\tau$. \end{assumption}

Overall, Assumption \ref{assump: A2} requires that the FIE has a
property that mimics the i-IOSS property of the system, while Assumption
\ref{assump: A3} ensures that the FIE is more sensitive than the
system to the uncertainties so that accurate inference of the state
is possible. The interpretation of the relatively more obscure Assumption
\ref{assump: A3} becomes clear with a concrete realization below. 
\begin{lem}
\label{lem: Special-form-of-A3}(Concrete realization of Assumption
\ref{assump: A3}) Assumption \ref{assump: A3} is true if the $\KL$
function $\underbar{\ensuremath{\rho}}$ satisfies 
\begin{equation}
\underbar{\ensuremath{\rho}}(|\tilde{\pi}_{i}|,\,\tau)\ge\alpha(2|\tilde{\pi}_{i}|,\,\tau),\label{eq: key-stability-condition-variant}
\end{equation}
for all $\tilde{\boldsymbol{\pi}}_{0:2t}\in\Pi_{\delta_{0}}$ and
$\tau\in\mathbb{I}_{0:t}$. 
\end{lem}
\begin{IEEEproof}
It suffices to show that inequality \eqref{eq: Assumption 3} is satisfied.
Subject to \eqref{eq: key-stability-condition-variant}, we have 
\begin{align*}
 & \alpha\left(2\underbar{\ensuremath{\rho}}^{-1}\left(\rho(s,\tau),\tau'\right),\tau'\right)\le\underbar{\ensuremath{\rho}}\left(\underbar{\ensuremath{\rho}}^{-1}\left(\rho(s,\tau),\tau'\right),\tau'\right)=\rho(s,\tau).
\end{align*}
This implies that inequality \eqref{eq: Assumption 3} is satisfied
with $\bar{\alpha}:=\rho$, and hence completes the proof. 
\end{IEEEproof}
Robust stability of FIE/MHE can then be established under Assumptions
\ref{assump: A2} and \ref{assump: A3}. 
\begin{thm}
\label{thm:RGAS-of-FIE}The following two conclusions hold:

a) (RAS/RGAS/RES/RGES of FIE) The FIE is RAS if the system is i-IOSS and the cost function of
FIE satisfies Assumptions \ref{assump: A2} and \ref{assump: A3},
and is RGAS if further the two assumptions are valid for all $\tilde{\boldsymbol{\pi}}_{0:2t}\in\Pi$.
And if, furthermore, the $\KL$ functions in Assumption \ref{assump: A3}
admit specific forms as $\alpha(s,\tau)=csb^{\tau}$ and $\bar{\alpha}(s,\tau)=\bar{c}s\bar{b}^{\tau}$
with certain $c,\bar{c}>0$, $b,\bar{b}\in(0,1)$ and for all $\tau\ge0$
and $s$ in the applicable domain, then the FIE which is RAS or RGAS
will be RES or RGES, respectively. 

b) (pRAS/pRGAS/RES/RGES/RGAS of MHE) In the four cases of a), the corresponding MHE under additional
Assumption \ref{assump: A1} with a sufficiently long horizon is pRAS,
pRGAS, RES and RGES, respectively. And the pRGAS of MHE strengthens
to RGAS if the $\K$ functions $\alpha(\cdot,0)$ and $\bar{\alpha}(\cdot,0)$
in Assumption \ref{assump: A3} are globally Lipschitz continuous
at the origin.
\end{thm}
\begin{IEEEproof}
\textit{a) RAS/RGAS of FIE}. The global optimal solution of $\tilde{\boldsymbol{\pi}}_{0:2t}$
for the FIE is denoted as $\hat{\tilde{\boldsymbol{\pi}}}_{0:2t}$
(cf. Eq. \eqref{eq: pi's}), yielding a minimum cost $V_{t}^{o}$.
It follows that for all $t\ge0$, 
\begin{align*}
V_{t}^{o} & =V_{t}(\hat{\tilde{\boldsymbol{\pi}}}_{0:2t})\le V_{t}(\boldsymbol{\pi}_{0:2t})\\
 & \stackrel{\text{Assump. \ref{assump: A2}}}{\le}\max_{i\in\mathbb{I}_{0:2t}}\ensuremath{\rho}(\left|\pi_{i}\right|,\,t-\iota(\pi_{i})-1)=:\bar{V}_{t}
\end{align*}
Consequently, by Assumption \ref{assump: A2} we have $\underbar{\ensuremath{\rho}}(|\hat{\tilde{\pi}}_{i}|,\,t-\iota(\hat{\tilde{\pi}}_{i})-1)\le V_{t}^{o}\le\bar{V}_{t}$,
and further $|\hat{\tilde{\pi}}_{i}|\le\underbar{\ensuremath{\rho}}_{t}^{-1}(\bar{V}_{t},\,t-\iota(\hat{\tilde{\pi}}_{i})-1)$.

Since $\hat{\pi}_{i}=\pi_{i}-\hat{\tilde{\pi}}_{i}$ for each $i\in\mathbb{I}_{0:2t}$,
by applying the triangle inequality this implies that 
\begin{align}
\left|\hat{\pi}_{i}\right| & \le\left|\pi_{i}\right|+\left|\hat{\tilde{\pi}}_{i}\right|\le\left|\pi_{i}\right|+\underbar{\ensuremath{\rho}}_{t}^{-1}(\bar{V}_{t},\,t-\iota(\hat{\tilde{\pi}}_{i})-1)\nonumber \\
 & \le2\left|\pi_{i}\right|\oplus2\underbar{\ensuremath{\rho}}_{t}^{-1}(\bar{V}_{t},\,t-\iota(\hat{\tilde{\pi}}_{i})-1),\label{eq:x0-x0^}
\end{align}
for all $i\in\mathbb{I}_{0:2t}$. Substitute \eqref{eq:x0-x0^} into
the i-IOSS property of \eqref{eq: concise definition - i-IOSS}, yielding
{\small{}{} 
\begin{align*}
 & \left|x(t;\,x_{0},\,\boldsymbol{w}_{0:t-1})-x(t;\,\hat{x}_{0},\,\hat{\boldsymbol{w}}_{0:t-1})\right|\\
 & \le\max_{i\in\mathbb{I}_{0:2t}}\alpha(\left|\hat{\pi}_{i}\right|,\,t-\iota(\hat{\pi}_{i})-1)
  \le\max_{i\in\mathbb{I}_{0:2t}}\alpha\left(2\left|\pi_{i}\right|,\,t-\iota(\hat{\pi}_{i})-1\right)\\
 & \quad\quad\quad\quad\quad\oplus\max_{i\in\mathbb{I}_{0:2t}}\alpha\left(2\underbar{\ensuremath{\rho}}^{-1}(\bar{V}_{t},\,t-\iota(\hat{\tilde{\pi}}_{i})-1),\,t-\iota(\hat{\pi}_{i})-1\right)\\
 & =\max_{i\in\mathbb{I}_{0:2t}}\alpha\left(2\left|\pi_{i}\right|,\,t-\iota(\hat{\pi}_{i})-1\right)\\
 & \quad\oplus\max_{i\in\mathbb{I}_{0:2t}}\max_{j\in\mathbb{I}_{0:2t}}\alpha\left(\begin{array}{c}
2\underbar{\ensuremath{\rho}}^{-1}\left(\begin{array}{c}
\ensuremath{\rho}(\left|\pi_{j}\right|,\,t-\iota(\pi_{j})-1),\\
t-\iota(\hat{\tilde{\pi}}_{i})-1
\end{array}\right),\\
t-\iota(\hat{\pi}_{i})-1
\end{array}\right)\\
 & \le\max_{i\in\mathbb{I}_{0:2t}}\left(\alpha\left(2\left|\pi_{i}\right|,\,t-\iota(\pi_{i})-1\right)\oplus\bar{\alpha}\left(\left|\pi_{i}\right|,\,t-\iota(\pi_{i})-1\right)\right)\\
 & =:\max_{i\in\mathbb{I}_{0:2t}}\beta\left(\left|\pi_{i}\right|,\,t-\iota(\pi_{i})-1\right),
\end{align*}
}where the equality $\iota(\hat{\pi}_{\cdot})=\iota(\pi_{\cdot})=\iota(\hat{\tilde{\pi}}_{\cdot})$
and Assumption \ref{assump: A3} have been used to derive the last
inequality, and $\beta\left(|\pi_{i}|,\tau\right):=\alpha\left(2|\pi_{i}|,\tau\right)\oplus\bar{\alpha}\left(|\pi_{i}|,\tau\right)$
for all $\boldsymbol{\pi}_{0:2t}\in\Pi_{\delta_{0}}$ and $\tau\ge0$.
Since $\beta$ is a $\KL$ function, the FIE is RAS by definition.
With the same reasoning, the conclusion immediately extends to that
the FIE is RGAS when Assumptions \ref{assump: A2} and \ref{assump: A3}
are applicable to the entire problem domain as specified by $\tilde{\boldsymbol{\pi}}_{0:2t}\in\Pi$.

\textit{RES/RGES of FIE}. If in addition, $\alpha\left(2s,\,\tau\right)=2csb^{\tau}$
and $\bar{\alpha}\left(s,\,\tau\right)=\bar{c}s\bar{b}^{\tau}$ with
certain $c,\bar{c}>0$ and $b,\bar{b}\in(0,1)$ and for all $\tau\ge0$
and $s$ in an applicable domain, then it is immediate that $\alpha\left(2s,\,\tau\right)\oplus\bar{\alpha}\left(s,\,\tau\right)\le\bar{c}'s\bar{b}'^{\tau}=:\beta(s,\,\tau)$,
with $\bar{c}':=\max\{2c,\,\bar{c}\}$ and $\bar{b}':=\max\{b,\,\bar{b}\}$.
Here $\beta$ is a $\KL$ function in an exponential form, so the
FIE which is RAS (or RGAS) in this case will be RES (or RGES) by Lemma
\ref{def: RGAS-1}. 

\textit{b) pRAS/pRGAS/RES/RGES/RGAS of MHE}. The conclusions are direct
application of RAS/RGAS/RES/RGES of FIE and Lemma \ref{lem: FIE-MHE-aRGAS and RGES}
and \ref{lem: FIE-MHE-RAS-aRAS}, subject to the given conditions.
\end{IEEEproof}
Next, we present a lemma indicating that Assumptions \ref{assump: A2}
and \ref{assump: A3} do not impose special difficulty as there always
exists a cost function satisfying both of them if the system is i-IOSS. 
\begin{lem}
\label{lem: Satisfaction-of-A2A3}(Satisfaction of Assumptions \ref{assump: A2}
and \ref{assump: A3}) If the system is i-IOSS as per \eqref{eq: concise definition - i-IOSS},
then, given any $t\ge0$ and $\underbar{\ensuremath{\rho}}\in\KL$
satisfying $\underbar{\ensuremath{\rho}}(|\tilde{\pi}_{i}|,\,\tau)\ge\alpha(2|\tilde{\pi}_{i}|,\,\tau)$
for all $\tilde{\boldsymbol{\pi}}_{0:2t}\in\Pi_{\delta_{0}}$ and
$\tau\in\mathbb{I}_{0:t}$, it is feasible to specify the cost function
of FIE as 
\begin{equation}
V_{t}(\tilde{\boldsymbol{\pi}}_{0:2t}):=\max_{i\in\mathbb{I}_{0:2t}}\underbar{\ensuremath{\rho}}(\left|\tilde{\pi}_{i}\right|,\,t-\iota(\tilde{\pi}_{i})-1),\label{eq: feasible form of FIE cost func}
\end{equation}
such that Assumptions \ref{assump: A2} and \ref{assump: A3} hold
true. The same form of cost function remains valid for Assumptions
\ref{assump: A2} and \ref{assump: A3} to hold globally with $\tilde{\boldsymbol{\pi}}_{0:2t}\in\Pi$.
In either the local or the global case, if the system is exp-i-IOSS,
then the $\KL$ function $\underbar{\ensuremath{\rho}}$ admits a
form as $\underbar{\ensuremath{\rho}}(s,\,\tau)=c|\tilde{\pi}_{i}|b^{\tau}$
with certain $c>0$ and $b\in(0,1)$ and for $\tau\ge0$ and $\tilde{\pi}_{i}$
in the applicable domain. 
\end{lem}
\begin{IEEEproof}
Given the cost function specified as per \eqref{eq: feasible form of FIE cost func},
Assumption \ref{assump: A2} is automatically met. With $\underbar{\ensuremath{\rho}}(|\tilde{\pi}_{i}|,\,\tau)\ge\alpha(2|\tilde{\pi}_{i}|,\,\tau)$
for all $\tilde{\boldsymbol{\pi}}_{0:2t}\in\Pi_{\delta_{0}}$ and
$\tau\in\mathbb{I}_{0:t}$, Assumption \ref{assump: A3} is also met
by Lemma \ref{lem: Special-form-of-A3}. While, existence of such
$\KL$ function $\underbar{\ensuremath{\rho}}$ is guaranteed as it
is always feasible to let $\underbar{\ensuremath{\rho}}(|\tilde{\pi}_{i}|,\,\tau):=\alpha(2|\tilde{\pi}_{i}|,\,\tau)$.
Since the reasoning approach does not rely on the domain of $|\tilde{\pi}_{i}|$,
the conclusion remains valid if Assumptions \ref{assump: A2} and
\ref{assump: A3} are extended to the entire problem domain with $\tilde{\boldsymbol{\pi}}_{0:2t}\in\Pi$.

If the system is exp-i-IOSS with $\alpha(|\tilde{\pi}_{i}|,\,\tau):=c'|\tilde{\pi}_{i}|b^{\tau}$
for certain $c'>0$ and $b\in(0,1)$, then it is valid to let $\underbar{\ensuremath{\rho}}(|\tilde{\pi}_{i}|,\,\tau):=\alpha(2|\tilde{\pi}_{i}|,\,\tau)=c|\tilde{\pi}_{i}|b^{\tau}$
with $c:=2c'$ for all $\tau\ge0$ and $\tilde{\pi}_{i}$ in the applicable
local/global domain. 
\end{IEEEproof}
Lemma \ref{lem: Satisfaction-of-A2A3} indicates that a valid cost
function can always be designed from the i-IOSS bound function $\alpha$,
and hence implies an important result below. 
\begin{cor}
\label{cor: if-and-only-if condition} There exists a cost function
for FIE to be RGAS (or RGES) if and only if the system is i-IOSS (or
exp-i-IOSS). 
\end{cor}
\begin{IEEEproof}
\textit{Sufficiency}. By Lemma \ref{lem: Satisfaction-of-A2A3}, the
system being i-IOSS implies that the FIE admits a cost function such
that Assumptions \ref{assump: A2} and \ref{assump: A3} hold true,
which consequently implies that the FIE is RGAS by Theorem \ref{thm:RGAS-of-FIE}.
When the system is exp-i-IOSS, the conclusion trivially strengthens
to that the FIE is RGES by following the same approach of reasoning.

\textit{Necessity}. The necessity in the RGAS case has been proved
in existing literature, e.g., \cite[Proposition 4.6]{rawlings2020model}
and \cite[Proposition 2.4]{allan2020lyapunov}, while the proof in
the RGES case is implied by the same proof there. 
\end{IEEEproof}
By Lemma \ref{lem: FIE-MHE-aRGAS and RGES}, it follows immediately
from Corollary \ref{cor: if-and-only-if condition} that there exists
a cost function for the MHE to be pRGAS (or RGES) if the system is
i-IOSS (or exp-i-IOSS). 
\begin{rem}
\label{rem: relation-to-most-recent-work2}The conclusion of Corollary
\ref{cor: if-and-only-if condition} coincides with a key finding
reported in a latest paper \cite{knuefer2021nonlinear}, which had
been submitted for review. The derivation approaches are, however,
quite different. Here, we apply the reasoning approach of \cite{Hu2017robust},
focusing on developing most general conditions for robust stability
of FIE, and the aforementioned conclusion appears as a corollary for
an endeavour to understand the developed conditions. In contrast,
reference \cite{knuefer2021nonlinear} reaches the conclusion by starting
with a particular cost function which is constructed directly from
the i-IOSS property of the system and is not necessarily the only
form admitted by our derived conditions. 
\end{rem}

\section{Conclusion\label{sec:Conclusion}}

This work proved that robust global (or local) asymptotic stability
of full information estimation (FIE) implies \textit{practical} robust
global (or local) asymptotic stability of moving horizon estimation
(MHE) which implements a sufficiently long horizon. The ``practical''
becomes exact if the FIE admits a certain Lipschitz continuity. In
both exact and inexact cases, explicit ways were also provided of
computing a sufficient horizon size for a robustly stable MHE. With
the revealed implication, sufficient conditions for the MHE to be
robustly stable were derived by firstly developing those for ensuring
robust stability of the corresponding FIE. A particular realization
of these conditions indicates that the system being i-IOSS is not
only necessary but also sufficient to ensure the existence of a robustly
globally asymptotically stable of FIE. With the revealed implication,
the sufficiency remains valid for existence of a practically robustly
stable MHE.

Since it is generic and relies only on robust stability of the corresponding
FIE, the revealed stability link implies that existing conditions
which ensure robust stability of FIE can all be inherited to establish
that of the corresponding MHE, but also paves the way for derivation
of new sufficient conditions via deeper analysis of FIE. This may also contribute to developing robustly stable MHE which
exploits sub-optimal solutions for resource-constrained or faster
estimations. Readers are referred to Sec. VII of \cite{Hu2017robust}
for related discussions, and \cite{wynn2014convergence,alessandri2017fast,gharbi2020anytime,schiller2020robust}
for some recent developments in this line of research.

\section*{Appendix A. Derivative of the Error Bound Factor w.r.t. the Moving
Horizon Size $T$}

The symbols mostly come from Lemma \ref{lem: FIE-MHE-RAS-aRAS}.a)
and its proof. Given $T\ge\underline{T}(\epsilon)$, it is easy to
show that the estimate error bound in \eqref{eq: error bound} can
be tightened by replacing $\eta^{\left\lfloor \frac{t}{T}\right\rfloor }$
with $\left(\frac{\eta\beta_{x}(s,T)}{\beta_{x}(s,\underline{T})}\right)^{\left\lfloor \frac{t}{T}\right\rfloor }$.
This appendix analyzes the derivative of this new bound factor w.r.t.
$T$, in order to understand monotonicity of the bound w.r.t. the
horizon size $T$.

Firstly, we show that there exists $\phi\in\LL$ such that $\phi(T-\underline{T})=\frac{\beta_{x}(s,T)}{\beta_{x}(s,\underline{T})}$
for all $T\ge\underline{T}$ . By Lemma 8 in \cite{sontag1998comments},
given $\beta_{x}\in\KL$, there exists $\alpha\in\K$ and $\varphi\in\LL$
such that $\beta_{x}(s,t)\le\alpha(s)\varphi(t)$, for all $s,t\ge0$.
This implies that the RGAS property of FIE can equivalently be expressed
with the product of $\K$ and $\LL$ functions. Therefore, it does
not lose generality by assuming $\beta_{x}(s,t):=\alpha(s)\varphi(t)$.
Consequently, $\frac{\beta_{x}(s,T)}{\beta_{x}(s,\underline{T})}=\frac{\varphi(T)}{\varphi(\underline{T})}$.
Since $\varphi$ is a $\LL$ function, given $\underline{T}\ge0$,
it is feasible to define $\phi(T-\underline{T})=\frac{\varphi(T)}{\varphi(\underline{T})}$
and $\phi$ will be a $\LL$ function for all $T\ge\underline{T}$.

Next, we derive the derivative of $\left(\eta\phi(T-\underline{T})\right)^{\left\lfloor \frac{t}{T}\right\rfloor }$
w.r.t. $T$. By definition we have $0<\phi(T-\underline{T})\le1$
for all $T\ge\underline{T}$. Let $z:=\ln\left(\eta\phi(T-\underline{T})\right)^{\left\lfloor \frac{t}{T}\right\rfloor }=\left\lfloor \frac{t}{T}\right\rfloor \ln(\eta\phi(T-\underline{T}))$.
Hence $\left(\eta\phi(T-\underline{T})\right)^{\left\lfloor \frac{t}{T}\right\rfloor }=e^{z}$.
The chain rule of derivatives implies

{\small{}{} 
\begin{align*}
 & \frac{d\left(\eta\phi(T-\underline{T})\right)^{\left\lfloor \frac{t}{T}\right\rfloor }}{dT}=e^{z}\frac{dz}{dT}\\
 & =\left(\eta\phi(T-\underline{T})\right)^{\left\lfloor \frac{t}{T}\right\rfloor }\cdot\left(\begin{array}{c}
\frac{d\left\lfloor \frac{t}{T}\right\rfloor }{dT}\ln\eta\phi(T-\underline{T})\\
+\frac{\left\lfloor \frac{t}{T}\right\rfloor }{\phi(T-\underline{T})}\frac{d\phi(T-\underline{T})}{dT}
\end{array}\right)\\
 & =\frac{\left\lfloor \frac{t}{T}\right\rfloor \left(\eta\phi(T-\underline{T})\right)^{\left\lfloor \frac{t}{T}\right\rfloor }}{\phi(T-\underline{T})}\cdot\left(\begin{array}{c}
\frac{\phi(T-\underline{T})}{\left\lfloor \frac{t}{T}\right\rfloor }\frac{d\left\lfloor \frac{t}{T}\right\rfloor }{dT}\ln\eta\phi(T-\underline{T})\\
+\frac{d\phi(T-\underline{T})}{dT}
\end{array}\right).
\end{align*}
}The first factor of the last expression is positive, so the sign
of the derivative is uniquely determined by the second factor. When
$t$ is much larger than $T$, we have $\left\lfloor \frac{t}{T}\right\rfloor \rightarrow\frac{t}{T}$.
Then, the second factor approximates to 
\begin{align*}
 & \frac{d\phi(T-\underline{T})}{dT}-\frac{1}{T}\phi(T-\underline{T})\ln\eta\phi(T-\underline{T})\\
 & =\frac{1}{\varphi(\underline{T})}\left(\frac{d\varphi(T)}{dT}-\frac{\varphi(T)}{T}\ln\frac{\eta\varphi(T)}{\varphi(\underline{T})}\right)
\end{align*}
As $\frac{1}{\varphi(\underline{T})}$ is positive, it is sufficient
to analyze the sign of 
\[
\kappa(T):=\frac{d\varphi(T)}{dT}-\frac{\varphi(T)}{T}\ln\frac{\eta\varphi(T)}{\varphi(\underline{T})}.
\]

As examples, two particular forms of $\varphi$ are considered. Firstly,
let us consider the exponential form: $\varphi(T):=b_{1}^{T}$ with
$b_{1}\in(0,\,1)$. Then, $\kappa(T)=\frac{b_{1}^{T}}{T}\ln\frac{b_{1}^{\underline{T}}}{\eta}$.
So, $\kappa(T)<0$ and hence $\left(\eta\phi(T-\underline{T})\right)^{\left\lfloor \frac{t}{T}\right\rfloor }$
will be strictly decreasing w.r.t. $T$ if $T\ge\underline{T}\ge1+\left\lceil \log_{b_{1}}\eta\right\rceil $.
Secondly, let us consider the fractional form: $\varphi(T):=(T+1)^{-b_{2}}$
with $b_{2}>0$. Then, $\kappa(T)=\frac{-b_{2}}{(T+1)^{b_{2}}}\left(\frac{1}{T+1}-\frac{1}{T}\ln\frac{\sqrt[b_{2}]{\eta}(T+1)}{\underline{T}+1}\right).$
So, $\kappa(T)<0$ and hence $\left(\eta\phi(T-\underline{T})\right)^{\left\lfloor \frac{t}{T}\right\rfloor }$
will again be strictly decreasing w.r.t. $T$ if $T\ge\underline{T}>\sqrt[b_{2}]{\eta}(T+1)e^{-\frac{T}{T+1}}-1$.
Since the condition is more restrictive in this case, the monotonic
decrement property seems harder to establish than that in the previous
case.

To conclude, the analysis shows that if the moving horizon size $T$
satisfies certain conditions (being large enough in general), the
bound of the estimate error as controlled by the factor $\left(\eta\phi(T-\underline{T})\right)^{\left\lfloor \frac{t}{T}\right\rfloor }$
will eventually decrease with $T$.

\bibliographystyle{IEEEtran}
\bibliography{Ref.bib}
%\bibliography{/Users/hwh/Dropbox/Research/Reference_database/Ref}
\end{document}